\documentclass{elegantpaper}

\usepackage{amsmath}
\usepackage{amsfonts}
\usepackage{amssymb}
\usepackage{mathrsfs}
\usepackage{graphicx}
\usepackage{xcolor}

\newcommand  \bglb {\big (}
\newcommand  \bgrb {\big )}

\newcommand   \scalars  {\mathbb F}

\newcommand  \dd   {\mathbf d}
\newcommand  \ee   {\mathbf e}

\newcommand  \xx   {\mathbf x}
\newcommand  \yy   {\mathbf y}
\newcommand  \zz   {\mathbf z}

\newcommand \norm[2] {\left \Vert #1 \right \basel{\Vert}{#2}}

\newcommand  \realNums       {\mathbb R}
\newcommand  \complexNums     {\mathbb C}

\newcommand   \naturalNums   {\mathbb N}
\newcommand   \indSetII      {\mathbb I}
\newcommand   \indSetJJ      {\mathbb J}

\newcommand  \basel[2]{#1_{_{#2}}}

\newcommand  \tab  {\hspace*{0.5cm}}

\newtheorem{Corollary}{\bf Corollary}[theorem]

\definecolor{titlecolor}{RGB}{144,48,48}
\definecolor{authorname}{RGB}{16,96,16}%
\definecolor{addrscolor}{RGB}{60,113,183}
\definecolor{secheader}{RGB}{16,32,128}%
\definecolor{refscolor}{RGB}{16,64,128}%

\begin{document}

\title{\textcolor{titlecolor}{\bf{Reflexivity of a Banach Space with a Countable Vector Space Basis}}}

\author{\\
\\
{\textcolor{authorname}{\Large{\bf{Michael Oser Rabin$^{1}$}}}} \tab
 {\textcolor{authorname}{and}} \tab
{\textcolor{authorname}{\Large{\bf{Duggirala Ravi$^{2}$}}}}\\
\\
 {\textcolor{addrscolor}{\scriptsize{\bf{$^{1}$E-mail ~: \tab rabin@math.seas.harvard.edu}}}} \\
 {\textcolor{addrscolor}{\scriptsize{\bf{$^{2}$E-mail ~: \tab ravi@gvpce.ac.in; ~~ drdravi2000@yahoo.com;
 ~~	 duggirala.ravi@rediffmail.com; ~~duggirala.ravi@yahoo.com} }}}
\\
\\
}

\date{}

\maketitle

\textcolor{secheader}{
\begin{abstract}
All most all the function spaces over real or complex domains and spaces of sequences, that arise in practice as examples of normed complete linear spaces (Banach spaces), are reflexive. These Banach spaces are dual to their respective spaces of continuous linear functionals over the corresponding Banach spaces. For each of these Banach spaces, a countable vector space basis exists, which is responsible for their reflexivity. In this paper, a specific criterion for reflexivity 
of a Banach space with a countable vector space basis is presented.
\end{abstract}
}

\begin{small}
\textcolor{secheader}{
}
\end{small}

\textcolor{titlecolor}{\section{\bf{Introduction}}}
All most all the function spaces over real or complex domains and spaces of sequences, that arise in practice as examples of normed complete linear spaces (Banach spaces), are reflexive. The topological property of being reflexive is that these Banach spaces are dual to their respective spaces of continuous linear functionals over the corresponding Banach spaces. For a reflexive Banach space, the double dual space is isometrically isomorphic to the Banach space, and they are essentially the same spaces. For each of these Banach spaces, a countable vector space basis exists, which is responsible for their reflexivity. A vector space basis for a Banach space can be easily formulated, taking the analogy from a Riesz basis for a Hilbert space. However, these Banach spaces are not plainly $L^{p}$-function or $\ell^{p}$-sequence spaces. For a Banach space with a countable vector space basis, the projection maps onto finite dimensional component subspaces are continuous and surjective, and hence open maps. The coefficient sequence of any vector in the Banach space, with respect to a basis, for which the dual basis functionals are normalized, can be shown to be bounded with respect to sup-norm, and the linear transformation mapping a vector to its coefficient sequence becomes continuous. The dual basis linear functionals become continuous and form a basis for the dual space.  A specific criterion for reflexivity of a Banach space with a countable vector space basis is presented.

\textcolor{titlecolor}{\section{\label{Sec-Main}\bf{Main Results}}}

Let $\naturalNums$ be the set of positive integers, and let $\realNums$ and $\complexNums$ be the fields of real and  
complex numbers, respectively, equipped with absolute value norm, denoted by $\vert \cdot \vert$.
Let $\scalars$ be $\realNums$ or $\complexNums$.

 A a topologically complete normed linear (vector) space over $\scalars$ is called a {\em Banach space}. 
Let $\ell^{p}\bglb\naturalNums, \scalars\bgrb $ be the Banach space of $\scalars$-valued sequences,
equipped with the $p$-norm, for $1 \leq p \leq \infty$.  Let $\basel{\ee}{i}$ be the standard
Euclidean vector, with $j$-th component of $\basel{\ee}{i}$ equal to the Kronecker $\basel{\delta}{i,\, j}$,
for $i, j \in \naturalNums$. The set $\{\basel{\ee}{i}\, : \,  i \in \naturalNums\}$ forms a {\em basis}
for $\ell^{p}\bglb\naturalNums, \scalars\bgrb $, for $1 \leq p \leq \infty$, and the basis is countable.

\textcolor{titlecolor}{\subsection{\label{Sec-Countable-Basis-Banach-Spaces}Countable Basis Banach Spaces}}

Let ${\mathcal B}$ be a Banach space, with norm 
$\norm{\cdot}{\mathcal B}$, and let 
$\{\basel{\xi}{i}\, :\, i \in \naturalNums\}\subset {\mathcal B }$. 
If every element $\xx \in {\mathcal B}$ can be expressed as a linear combination
$\xx = \sum_{i \in \naturalNums} \basel{c}{i} \basel{\xi}{i}$, for some scalars
$\basel{c}{i} \in \scalars$, where $i \in \naturalNums$, then the Banach space
is said to admit a {\em countable basis}, which is
 $\{\basel{\xi}{i}\, :\, i \in \naturalNums\}$.
 By the uniqueness of the linear combination, it is assumed that if 
$\sum_{i \in \naturalNums} \basel{d}{i} \basel{\xi}{i} = 0$,
for some scalars $\basel{d}{i} \in \scalars$,
then $\basel{d}{i} = 0$, for every $i \in \naturalNums$.
Componentwise addition and subtraction is required to hold.
The space of continuous linear functionals,
denoted by $\hat{\mathcal B}$,
is called the {\em dual} of ${\mathcal B}$,
with its dual norm $\norm{\cdot}{\hat{\mathcal B}}$.
The double dual, or {\em bidual}, is denoted by $\hat{\hat{\mathcal B}}$,
with its bidual norm $\norm{\cdot}{\hat{\hat{\mathcal B}}}$.
\\

\begin{proposition}\label{Prop-2.1}
Let ${\mathcal B}$ be a Banach space, with norm 
$\norm{\cdot}{\mathcal B}$, and let 
$\{\basel{\xi}{i}\, :\, i \in \naturalNums\}\subset {\mathcal B }$
be a basis for ${\mathcal B}$.
Let $\indSetII \subset \naturalNums$ and 
$\indSetJJ = \naturalNums ~ \backslash ~ \indSetII$,
and let ${\mathcal L}(\indSetII)$ and  ${\mathcal L}(\indSetJJ)$ be the
subspaces spanned by $\{\basel{\xi}{i}\, :\, i \in \indSetII\}$ and
$\{\basel{\xi}{j}\, :\, j \in \indSetJJ\}$, respectively. The following
statements hold:
\begin{enumerate}
\item Both ${\mathcal L}(\indSetII)$ and  ${\mathcal L}(\indSetJJ)$
      are closed linear subspaces of ${\mathcal B }$, with respect to
      $\norm{\cdot}{\mathcal B}$.
\item If at least one of the index sets $\indSetII $ and $\indSetJJ $ is finite,
      then ${\mathcal B} = {\mathcal L}(\indSetII) \oplus {\mathcal L}(\indSetJJ)$,
      {\em {i.e.}}, every vector $\xx \in {\mathcal B}$ can be expressed as the
      sum of vectors $\xx = \yy + \zz$, for some $\yy \in {\mathcal L}(\indSetII)$
       and  $\zz \in {\mathcal L}(\indSetJJ)$, uniquely.
\item If at least one of the index sets $\indSetII $ and $\indSetJJ $ is finite,
      the projection operators $\basel{\mathsf P}{\indSetII}$ and 
      $\basel{\mathsf P}{\indSetJJ}$,  defined by
      $\basel{\mathsf P}{\indSetII}(\xx) = \yy$ and
      $\basel{\mathsf P}{\indSetJJ}(\xx) = \zz$,
      for $\xx \in {\mathcal B}$, with $\yy \in {\mathcal L}(\indSetII)$
      and  $\zz \in {\mathcal L}(\indSetJJ)$, such that  $\xx = \yy + \zz$,
      are both well defined and continuous.
             
\end{enumerate}
\end{proposition}
\proof If $\indSetII = \emptyset$ or $\indSetJJ = \emptyset$,
the statements vacuously hold. 

For a singleton set $\{ i \}$, for some $i \in \naturalNums$,
${\mathcal L}(\{i\}) = \{ c\basel{\xi}{i} \,:\,  c \in  \scalars\}$,
and $\norm{c\basel{\xi}{i}}{\mathcal B} = \vert c \vert \cdot \norm{\basel{\xi}{i}}{\mathcal B}$,
whence ${\mathcal L}(\{i\})$ is topologically closed in ${\mathcal B}$
with respect to $\norm{\cdot}{\mathcal B}$.
If $\xx = \sum_{j \in \naturalNums} \basel{c}{j}\basel{\xi}{j} \in {\mathcal B}$,
for some scalars $\basel{c}{j} \in \naturalNums$, then
$\yy = \basel{c}{i} \basel{\xi}{i} \in {\mathcal L}(\{i\})$,
$\zz = \sum_{j \in \naturalNums\, \backslash \, \{i\}} \basel{c}{j}\basel{\xi}{j}
 \in {\mathcal L}(\naturalNums\, \backslash \, \{i\})$, and the expression
$\xx = \yy + \zz$ is uniquely determined. Therefore,
${\mathcal L}(\naturalNums\, \backslash \, \{i\})$
is a closed linear subspace of ${\mathcal B}$,
with respect to $\norm{\cdot}{\mathcal B}$, as well.

Now ${\mathcal L}(\indSetII) = \bigcap_{j \in \indSetJJ} 
{\mathcal L}(\naturalNums\, \backslash \, \{j\})$
and ${\mathcal L}(\indSetJJ) = \bigcap_{i \in \indSetII} 
{\mathcal L}(\naturalNums\, \backslash \, \{i\})$,
hence both are closed.  The remaining part is obvious.                    \qed
\\

 For a Banach space  ${\mathcal B}$, with a countable basis 
$\{\basel{\xi}{i}\, :\, i \in \naturalNums\}\subset {\mathcal B }$,
the algebraic dual basis is the set of linear functionals 
$\{\basel{\hat{\xi}}{i}\, :\, i \in \naturalNums\}$,
defined by their action on the basis vectors by the 
condition that
$\basel{\hat{\xi}}{i}(\basel{\xi}{j}) = \basel{\delta}{i,\, j}$,
for $i, \, j \in \naturalNums$.
Proposition 1 (item 3) shows that $\basel{\hat{\xi}}{i}$
is continuous, hence  bounded, {\em {i.e.}}, 
$1 \leq \norm{\basel{\hat{\xi}}{i}}{\hat{\mathcal B}} < \infty$,
for every $i \in \naturalNums$.
Let $\basel{\eta}{i} = 
\norm{\basel{\hat{\xi}}{i}}{\hat{\mathcal B}}\basel{\xi}{i}$
and $\basel{\hat{\eta}}{i} = 
\frac{\basel{\hat{\xi}}{i}}{\norm{\basel{\hat{\xi}}{i}}{\hat{\mathcal B}}}$,
for $i \in \naturalNums$.
Then, $\norm{\basel{\hat{\eta}}{i}}{\hat{\mathcal B}} = 1$ and 
$\basel{\hat{\eta}}{i}(\basel{\eta}{j}) = \basel{\delta}{i,\, j}$\, ,
for $i,\, j \in \naturalNums$. Moreover,
$\{\basel{\hat{\eta}}{i}\, :\, i \in \naturalNums \} $
forms a basis for $\hat{\mathcal B}$ :   if
$\hat{f} \in \hat{\mathcal B}$, then,
by the algebraic action of $\hat{f}$ on the
vector space basis
$\{\basel{\eta}{i}\, :\, i \in \naturalNums\}$
for ${\mathcal B}$,  it holds that
$\hat{f}= \sum_{i \in \naturalNums} \hat{f}(\basel{\eta}{i}) \basel{\hat{\eta}}{i}$,
and the uniqueness of the algebraic  expression is obvious.
\\

\begin{proposition}
Let ${\mathcal B}$ be a Banach space, with a countable basis,
$\{\basel{\eta}{i}\, :\, i \in \naturalNums\}$, such that
$\{\basel{\hat{\eta}}{i}\, :\, i \in \naturalNums\}$ is
a normalized basis for $\hat{\mathcal B}$, {\em i.e.},
 $\norm{\basel{\hat{\eta}}{i}}{\hat{\mathcal B}} = 1$
 and $\{\basel{\hat{\eta}}{i}\, :\, i \in \naturalNums \} $
forms a basis for $\hat{\mathcal B}$, and let
 $\xx = \sum_{i \in \naturalNums} \basel{c}{i} \basel{\eta}{i}$,
for some $\basel{c}{i} \in \scalars$, where $i\in \naturalNums$.
Then,
$\sup_{i \in \naturalNums} \vert \basel{c}{i} \vert \leq  \norm{\xx}{\mathcal B}$.
\end{proposition}

\proof  Let $\hat {T}$ be the formal linear mapping defined by
$\hat{T}(\dd) = \sum_{i \in \naturalNums} \basel{d}{i} \basel{\hat{\xi}}{i}$,
for $\dd =  (\basel{d}{1},\, \basel{d}{2}, \, \basel{d}{3},\,\ldots) \in 
\ell^{1}(\naturalNums,\, \scalars)$. For $m, \, n \in \naturalNums$,
and $\basel{g}{n} = \sum_{i = 1}^{n}\basel{d}{i}\basel{\hat{\xi}}{i}$,
the estimate $\norm{\basel{g}{m+n}-\basel{g}{n}}{\mathcal B}$
$ = $
$\norm{\sum_{i=n+1}^{m+n}\basel{d}{i}\basel{\hat{\eta}}{i}}{\mathcal B}$
$ \leq $
$\sup_{n+1 \leq i\leq m+n}\sum_{i=n+1}^{m+n}\vert \basel{d}{i} \vert $,
shows that $\{\basel{g}{i} \, :\, i \in \naturalNums\}$ is a Cauchy sequence
in $\hat{\mathcal B}$, and $\hat{T}(\dd)\in \hat{\mathcal B}$, with the norm of
$\hat{T}$, as  a  linear transformation from $\ell^{1}(\naturalNums,\, \scalars)$
into $\hat{\mathcal B}$, being at most $1$. Let
 $\hat{\hat{T}}\,:\,\hat{\hat{\mathcal B}}\,\to\, \ell^{\infty}(\naturalNums,\,\scalars)$
be the adjoint linear transformation of $\hat{T}$.
The linear transformation norm of $\hat{\hat{T}}$
coincides with that  of $\hat{T}$, and  $\hat{\hat{T}}$
becomes continuous. Let $\basel{\hat{\hat{\eta}}}{i}$ be the
natural embedding of $\basel{\eta}{i}$ in $\hat{\hat{\mathcal B}}$.
For $\xx = \sum_{i \in \naturalNums} \basel{c}{i} \basel{\eta}{i}$,
for some scalars $\basel{c}{i} \in  \scalars$, $i \in \naturalNums$,
the natural embedding of $\xx$ in $\hat{\hat{\mathcal B}}$ is
$\hat{\hat{\xx}} = \sum_{i \in \naturalNums} \basel{c}{i} \basel{\hat{\hat{\eta}}}{i}$,
which is an isometric isomorphism, by Hahn-Banach theorem,
 {\em {i.e.}},
  $\norm{\xx}{\mathcal B}$
$ =  $ 
$ \norm{\hat{\hat{x}}}{\hat{\hat{\mathcal B}}}$.
If $ \norm{\hat{\hat{x}}}{\hat{\hat{\mathcal B}}} \leq 1$,
then $\vert \hat{\hat{\xx}}\bglb \hat{T}(\basel{\ee}{i})\bgrb \vert$
$ = $
 $ \vert \hat{\hat{\xx}}(\basel{\hat{\eta}}{i}) \vert $
$ = \vert \basel{c}{i} \vert  \leq 1$, for every $i \in \naturalNums$.     
By the linearity of of the  natural embedding, the contention follows.          \qed
 \\

\textcolor{titlecolor}{\subsection{\label{Sec-Reflexivity}Reflexivity of a Countable Basis Banach Space}}
 
Let ${\mathcal B}$ be a Banach space with a countable basis
 $\{\basel{\xi}{i} \,:\, i\in \naturalNums\}$, and
  $\hat{\mathcal B}$ be the dual space of with ${\mathcal B}$
  generated by the dual basis
 $\{\basel{\hat{\xi}}{i} \,:\, i\in \naturalNums\}$.
 In this subsection, three specific assumptions regarding
 structure of ${\mathcal B}$ and $\hat{\mathcal B}$ are
 assumed to hold good, for the theorem of this section.
\\

\noindent {\bf{Assumption 1.}} \tab
 For every $\xx \in {\mathcal B}$, with
 $\xx = \sum_{i \in \naturalNums}  \basel{c}{i} \basel{\xi}{i}$,
 for some scalars $ \basel{c}{i} \in  \scalars$,
 if $\basel{\yy}{n} = \sum_{i = 1}^{n}  \basel{c}{i} \basel{\xi}{i}$,
 then  $\norm{\basel{\yy}{n}}{\mathcal B} \leq \norm{\xx}{\mathcal B}$,
 for every $n \in \naturalNums$.
 \\
 
 The assumption just stated can hold,  for a very large collection of Banach spaces.
 It is  particularly a nondecreasing norm, with addition of more terms :
 $\norm{\basel{\yy}{n}}{\mathcal B} \leq \norm{\basel{\yy}{n+1}}{\mathcal B}$,
 for each $n \in \naturalNums$.  In general, for arbitrary scalars
 $\basel{c}{i} \in \scalars$, with $\vert\basel{c}{i}\vert$ bounded
 and  $\basel{\yy}{n} = \sum_{j = 1}^{n}  \basel{c}{j} \basel{\xi}{j}$,
 where $i,\, n \in \naturalNums$, it may not be true that 
 $\{\norm{\basel{\yy}{n}}{\mathcal B}\, :\, n \in \naturalNums\}$
 is a bounded set of nonnegative real numbers, even though
$ \norm{\basel{\yy}{n}}{\mathcal B}$ is nondecreasing, for $n \in \naturalNums$.
 \\

\noindent {\bf{Assumption 2.}} \tab
 For any scalars $\basel{c}{i} \in \scalars$, where $i \in \naturalNums$, 
 with $\basel{\yy}{n} = \sum_{i = 1}^{n}  \basel{c}{i} \basel{\xi}{i}$,
 if $\norm{\basel{\yy}{n}}{\mathcal B} \leq b$,
 for some fixed $b > 0$ and every $n \in \naturalNums$,
 then $ \sum_{i \in \naturalNums}  \basel{c}{i} \basel{\xi}{i}\in {\mathcal B}$,
 and if $\xx =  \sum_{i \in \naturalNums}  \basel{c}{i} \basel{\xi}{i}
 \in {\mathcal B}$, then  $ \norm{\xx}{\mathcal B} \leq b$.
 \\
 
 The assumption just stated can hold,  for a very large collection of Banach spaces.
 By the embedding of 
 $\basel{\yy}{n} = \sum_{i = 1}^{n}  \basel{c}{i} \basel{\xi}{i}\in {\mathcal B}$
 into the double dual
  $\basel{\hat{\hat{\yy}}}{n} = \sum_{i = 1}^{n}  \basel{c}{i} \basel{\hat{\hat{\xi}}}{i}\in \hat{\hat{\mathcal B}}$,
  followed by an appeal to the compactness of the closed convex unit sphere of
  $\hat{\hat{\mathcal B}}$, centered at the origin, with respect to weak$^{*}$
  topology of  $\hat{\hat{\mathcal B}}$,
 if $\norm{\basel{\yy}{n}}{\mathcal B} \leq b$,
 then the sequence 
   $\basel{\hat{\hat{\yy}}}{n} = \sum_{i = 1}^{n}  \basel{c}{i} \basel{\hat{\hat{\xi}}}{i}$
   converges to
      $\hat{\hat{\zz}} = \sum_{i \in \naturalNums}  \basel{c}{i} \basel{\hat{\hat{\xi}}}{i}$
   with $\norm{\hat{\hat{\zz}}}{\hat{\hat{\mathcal B}}} \leq b$.
  The assumption implies that there exists $\xx \in {\mathcal B}$,
  such that the natural embedding of $\xx$ into $\hat{\hat{\mathcal B}}$
  is $\hat{\hat{\zz}}$.
 \\

\noindent {\bf{Assumption 3.}}  \tab  
 For every $\hat{f} \in \hat{\mathcal B}$, with
 $\hat{f} = \sum_{i \in \naturalNums}  \basel{d}{i} \basel{\hat{\xi}}{i}$,
 for some scalars $ \basel{d}{i} \in  \scalars$,
 it holds that
 \begin{equation}
 \lim_{n\rightarrow \infty}\norm{\sum_{i = n+1}^{\infty}
  \basel{d}{i} \basel{\hat{\xi}}{i}~~}{\hat{\mathcal B}} ~~  = ~~ 0
       \label{dual-space-approximable-by-finite-number-of-terms}
\end{equation}
  
  Let $\hat{\mathcal M} $
  $  =  $ 
  $\bigcup_{n \in \naturalNums} \hat{\mathcal L}(\{1,\, ... , \,n\})$,
 where $\hat{\mathcal L}(\{1,\, ... , \,n\})$ is the closed linear subspace
 generated by $\{\basel{\hat{\xi}}{i} \,:\, 1 \leq  i \leq  n\}$,
 for  $n \in \naturalNums$.  If the assumption just stated  holds,
 then $\hat{\mathcal B}$ is the topological closure of
 $\hat{\mathcal M} $, with respect to $\norm{\cdot}{\hat{\mathcal B}}$.
 The assumption is satisfied, if the coefficient sequence is in
 $\ell^{p}(\naturalNums,\,\scalars)$, for some finite $p \geq 1$,  possibly
 depending on $\hat{f}$, {\em{i.e.}}, $1 \leq  p < \infty$.
 The following is the main result.
 \\
 
 \begin{theorem}
 {\bf {(Rabin)}~~}
 Let ${\mathcal B}$ be a Banach space with a countable basis
 $\{\basel{\xi}{i} \,:\, i\in \naturalNums\}$.
 If the three assumptions  stated above hold, then  ${\mathcal B}$
 is reflexive, {\em {i.e.}}, the double dual $\hat{\hat{\mathcal B}}$
 is isometrically isomorphic to  ${\mathcal B}$.
 \end{theorem}
 \proof Let $S(r) =  \{\xx \in {\mathcal B}\,:\, \norm{\xx}{\mathcal B} < r\}$,
 for $r > 0$, be the open subset of  ${\mathcal B}$, of radius $r$,
 centered at the origin, and
 $\overline{S(r)} =  \{\xx \in {\mathcal B}\,:\, \norm{\xx}{\mathcal B} \leq r\}$,
  for $r > 0$, be the closed subset of  ${\mathcal B}$, of radius $r$,
  centered at the origin.

 Let $\basel{V}{i}$
$ = $
$ \{ \basel{d}{i} \in \scalars\,:\,\vert \basel{d}{i} \vert$
$ < $
$ \norm{\basel{\hat{\xi}}{i}}{\hat{\mathcal  B}}\}$
be the open set of $\scalars$,
of radius $\norm{\basel{\hat{\xi}}{i}}{\hat{\mathcal  B}}$,
centered at the origin, and
$\overline{\basel{V}{i}}$
$ = $
$ \{ \basel{d}{i} \in \scalars\,:\,\vert \basel{d}{i} \vert$
$ \leq $
$ \norm{\basel{\hat{\xi}}{i}}{\hat{\mathcal  B}}\}$
be the closed set of $\scalars$, 
of radius $\norm{\basel{\hat{\xi}}{i}}{\hat{\mathcal  B}}$,
centered at the origin.

 The product topology on $\overline{S(1)}$ is generated by the basic open sets
 \begin{scriptsize}
 \[
\left\{~~
\sum_{i \in \naturalNums}\basel{c}{i}\basel{\xi}{i} \in \overline{S(1)} ~~~~ :  ~~~~
\basel{c}{i} \in \basel{U}{i}\,, ~~  i \in \naturalNums ~~ \right \} 
\]
 \end{scriptsize}
where  $\basel{U}{i}$ is an open subset of $\overline{\basel{V}{i}}$,
  such that $\basel{U}{i} = \overline{\basel{V}{i}}$, for all but possibly
  finitely many indexes $i \in \naturalNums$.
  
Let $T\,:\,\overline{S(1)} ~ \to ~ \prod_{i\in\naturalNums}\overline{\basel{V}{i}}$ be the
on-to-one continuous function mapping
 $\xx = \sum_{i\in\naturalNums}\basel{c}{i}\basel{\xi}{i}$
to the coefficient sequence $(\basel{c}{1},\, \basel{c}{2},\, \ldots) $
$  \in  $
$\prod_{i\in\naturalNums}\overline{\basel{V}{i}}$.
Then, the basic open  set of $\overline{S(1)}$ as just defined 
is the set $T^{-1}\bglb \prod_{i \in \naturalNums}\basel{U}{i}\bgrb$, where
$\prod_{i \in \naturalNums}\basel{U}{i}$ is an open subset of the countable
infinite product $\prod_{i\in\naturalNums}\overline{\basel{V}{i}}$.

There are two parts in the proof: in the first part, $T\bglb~ \overline{S(1)} ~\bgrb$
is shown to be closed in $\prod_{i \in \naturalNums}\overline{\basel{V}{i}}$ (hence
$\overline{S(1)}$ becomes compact with respect to the product topology),
and in the second part, the weak topology of $\overline{S(1)}$ is shown
to  coincide with the product topology of $\overline{S(1)}$.
\\

\noindent{\sf{Part 1. ~~}}  Since the projections onto finite dimensional subspaces are continuous
   and also open maps, the set
 \begin{scriptsize}
 \[
\basel{\Gamma}{n} ~~~~ = ~~~~
\left\{~~(\basel{c}{1},\, ...\, ,\, \basel{c}{n}) \in \scalars^{n} ~~~~ :  ~~~~
\norm{\sum_{i\in \naturalNums}\basel{c}{i}\basel{\xi}{i}~~}{\mathcal B} \leq 1~~,~~~~\textrm{for some scalars}~~\basel{c }{i} \in  \scalars ~~ \textrm{and} ~~ i \geq n+1 \right \}
\]
 \end{scriptsize}
for $n \in \naturalNums$, is a closed subset of $\prod_{i = 1}^{n}\overline{\basel{V}{i}}$.
Now, the set $\basel{A}{n}  = \basel{\Gamma}{n}\times \prod_{i = n+1}^{\infty}\overline{\basel{V}{i}}$
is a  closed subset of $\prod_{i \in \naturalNums}\overline{\basel{V}{i}} $,
with respect to the product topology, for every $n \in \naturalNums$.
If $\norm{\sum_{i\in\naturalNums}\basel{c}{i}\basel{\xi}{i}}{\mathcal B} \leq 1$,
then $\norm{\sum_{i=1}^{n}\basel{c}{i}\basel{\xi}{i}}{\mathcal B} \leq 1$,
by Assumption 1, and therefore,
$T\bglb~ \overline{S(1)} ~\bgrb \subseteq \basel{A}{n}$, for every $n \in \naturalNums$.
Conversely, if
 $(\basel{c}{1},\, ...\, ,\, \basel{c}{n},\,\basel{c}{n+1},\, \ldots) \in \basel{A}{n}$,
 for every $n \in \naturalNums$, then   
 $\norm{\sum_{i=1}^{n}\basel{c}{i}\basel{\xi}{i}~~}{\mathcal B} \leq 1$,
  for every $n \in \naturalNums$,  by Assumption 1, and
   $\sum_{i\in \naturalNums}\basel{c}{i}\basel{\xi}{i}\in {\mathcal B}$,
   with  $\norm{\sum_{i\in \naturalNums}\basel{c}{i}\basel{\xi}{i}~~}{\mathcal B} \leq 1$,
    by Assumption 2.   Thus, 
 \begin{scriptsize}
 \[
 \bigcap_{n \in \naturalNums} \basel{A}{n}  ~~~~  =   ~~~~ T\bglb~ \overline{S(1)} ~\bgrb
 \]
 \end{scriptsize}
 and $T\bglb~ \overline{S(1)} ~\bgrb$ is a closed subset of 
 $\prod_{i \in \naturalNums}\overline{\basel{V}{i}}$, hence
 compact with respect to the product topology.        
 \\

\noindent{\sf{Part 2. ~~}} The weak topology on ${\mathcal B}$
 is generated by subbasic open sets 
 ${\mathcal V}(\hat{f},\, \zz, \, \rho)$, where
 $\hat{f}\in \hat{\mathcal B}$, $\zz \in {\mathcal B}$
 and $\rho > 0$, defined by the condition
 \[
 {\mathcal V}(\hat{f},\, \zz, \, \rho) ~~ = ~~
 \left\{
 \yy \in {\mathcal B}\,:\, \vert \hat{f}(\yy) -  \hat{f}(\zz)\vert < \rho
 \right\}
 \]
  For any $\xx \in  {\mathcal V}(\hat{f},\, \zz, \, \rho)$, 
   with $\delta = \delta(\hat{f}, \zz,\, \rho) =  \rho - \vert \hat{f}(\zz)-\hat{f}(\xx)\vert$,
   by triangle inequality,
  $ \vert \hat{f}(\yy) -  \hat{f}(\zz)\vert $
  $ \leq $
  $ \vert \hat{f}(\yy) -  \hat{f}(\xx)\vert $
  $ + $
  $ \vert \hat{f}(\xx) -  \hat{f}(\zz)\vert $,
  if $\yy \in  {\mathcal V}(\hat{f},\, \xx, \, \delta)$,
  then  $ \vert \hat{f}(\yy) -  \hat{f}(\zz)\vert $
      $ \leq $
      $\delta + \vert \hat{f}(\xx) -  \hat{f}(\zz)\vert  < \rho$,
  and  $\yy \in  {\mathcal V}(\hat{f},\, \zz, \, \rho)$, hence
 ${\mathcal V}(\hat{f},\, \xx, \, \delta)\subseteq  {\mathcal V}(\hat{f},\, \zz, \, \rho)$.
 
  Now, for the basic open set  ${\mathcal V}(\hat{f},\, \xx, \, \delta)$,
  let $\xx = \sum_{i\in \naturalNums}\basel{c}{i}\basel{\xi}{i}$ and
   $\hat{f} = \sum_{i\in \naturalNums}\basel{d}{i}\basel{\xi}{i}$,
   for some scalars $\basel{c}{i},\,\basel{d}{i} \in \scalars$,
   for $i\in \naturalNums$.    Since
   $\hat{f}(\xx) = \sum_{i\in \naturalNums}\basel{c}{i}\basel{d}{i}$,
   it follows that there exists $N(\hat{f},\,\delta)\in \naturalNums$,
   such that
    $   \norm{\sum_{i=n+1}^{\infty}\basel{d}{i}\basel{\hat{\xi}}{i}}{\hat{\mathcal B}}$
    $ < $
    $ \frac{\delta}{4}$,
   for every $n \geq N(\hat{f},\,\delta)$.  
   Let $\hat{g} = \sum_{i=1}^{n}\basel{d}{i}\basel{\hat{\xi}}{i}$,
   for some fixed $n \geq N(\hat{f},\,\delta)$.
   If $\vert\hat{g}(\yy) - \hat{g}(\xx)\vert < \frac{\delta}{4}$, then
   $\vert\hat{f}(\yy) - \hat{f}(\xx)\vert$
   $\leq$
      $\vert\hat{f}(\yy) - \hat{g}(\yy)\vert$
      $+$
    $\vert\hat{g}(\yy) - \hat{g}(\xx)\vert$
          $+$
    $\vert\hat{g}(\xx) - \hat{f}(\xx)\vert$
    $ <  $
    $\frac{\delta}{4}\bglb \norm{\xx}{\mathcal B}+\norm{\yy}{\mathcal B}\bgrb$
    $+$
    $\frac{\delta}{4}$
    $ <\delta$, whenever $\xx,\, \yy  \in \overline{S(1)}$.
    Thus, 
    ${\mathcal V}(\hat{g},\, \xx, \, \frac{\delta}{4}) \bigcap \overline{S(1)}$
    $\subseteq$
        ${\mathcal V}(\hat{f},\, \zz, \, \rho) \bigcap \overline{S(1)}$,
        for $\xx \in {\mathcal V}(\hat{f},\, \zz, \, \rho) \bigcap \overline{S(1)}$,
       with $\delta = \delta(\hat{f}, \zz,\, \rho) =  \rho - \vert \hat{f}(\zz)-\hat{f}(\xx)\vert$,
        and some appropriately chosen $\hat{g} \in \hat{\mathcal M}$.
   \\
    
   In the product topology, all the functionals in $\hat{\mathcal M}$
   remain  continuous, and the product topology on $\overline{S(1)}$
   coincides with the weak topology  on $\overline{S(1)}$.
   Thus,  $\overline{S(1)}$ is compact in the weak topology.       \qed
  \\

\paragraph{Observation}\tab The compactness of the set
$\prod_{i \in \naturalNums}\overline{\basel{V}{i}} $
can be proved without using the axiom of choice.
It is also possible to take ${\mathcal B}$ to be
the dual of  another Banach space ${\mathcal  N}$,
and apply the criterion to infer the reflexivity
of ${\mathcal N}$.

\textcolor{titlecolor}{
\section{\label{Sec-Misc}Miscellany and Applications} 
}
A vector space basis is very important for a vector space.  For an infinite
dimensional topological vector space, a precise formulation of a vector space
basis can be expressed, extending the definition of linear independence for
linear combinations of finitely many vectors to series. For a Banach space
with a countable vector space basis, some more interesting and useful facts
can be proved.
\\

\begin{theorem}
{\bf{(Weak Denseness of Finite Linear Combinations)}}\tab   
\label{Weak-Denseness-of-Finite-Linear-Combinations}
 Let ${\mathcal B}$ be a Banach space with a countable basis
 $\{\basel{\xi}{i} \,:\, i\in \naturalNums\}$, and
 ${\mathcal M} =  \bigcup_{n \in \naturalNums} {\mathcal L}(\{1,\, ... , \,n\})$,
 where ${\mathcal L}(\{1,\, ... , \,n\})$ is the closed linear subspace
 generated by $\{\basel{\xi}{i} \,:\, 1 \leq  i \leq  n\}$,
 for  $n \in \naturalNums$.
 Then, ${\mathcal M}$ is weakly dense in ${\mathcal B}$.
 \end{theorem}
 \proof The weak topology on ${\mathcal B}$ is generated by
 subbasic open sets  ${\mathcal V}(\hat{f},\, \xx, \, \epsilon)$,
 where $\hat{f}\in \hat{\mathcal B}$, $\xx \in {\mathcal B}$
 and $\epsilon > 0$, where 
 $ {\mathcal V}(\hat{f},\, \xx, \, \epsilon)$
 $ = $
$ \left\{
 \yy \in {\mathcal B}\,:\, \vert \hat{f}(\yy) -  \hat{f}(\xx)\vert < \epsilon
 \right\}$.

 Let $\xx = \sum_{i \in \naturalNums} \basel{c}{i}\basel{\xi}{i}\in {\mathcal B}$,
 for some scalars $\basel{c}{i} \in  \scalars$, for $i \in \naturalNums$. Let 
 $\basel{\yy}{n} = \sum_{i = 1}^{n} \basel{c}{i}\basel{\xi}{i}$, so that
 $\basel{\yy}{n} \in {\mathcal L}(\{1,\, ... , \,n\})$, for  $n \in \naturalNums$.
 For every $\hat{f} = \sum_{i \in \naturalNums} \basel{d}{i}\basel{\hat{\xi}}{i}$,
 for any scalars  $\basel{d}{i} \in  \scalars$, for $i \in \naturalNums$, and
 $\epsilon > 0$, there exists $N(\hat{f},\,\xx,\, \epsilon) \in \naturalNums$,
 such that $\vert \hat{f}(\basel{\yy}{n}) - \hat{f}(\xx)\vert $
 $=$
 $\vert\,    \sum_{i \geq n+1} \basel{c}{i} \basel{d}{i} \, \vert$
 $ < \epsilon$, for every $n \geq N(\hat{f},\,\xx,\, \epsilon)$.
 Now, for a basic open subset
 ${\mathcal V}(\basel{\hat{f}}{1},\, \xx, \, \basel{\epsilon}{1}) \cap \, ...\, 
\cap {\mathcal V}(\basel{\hat{f}}{m},\, \xx, \, \basel{\epsilon}{m})  $,
for some $m \in \naturalNums$,
 let $\delta = \min\{\basel{\epsilon}{i}\,:\, 1\leq i\leq m\}$
 and $\nu = \max \{N(\basel{\hat{f}}{i},\,\xx,\, \delta)\,:\, 1\leq i\leq   m\}$.
 Then,  $\vert \basel{\hat{f}}{i}(\basel{\yy}{n}) - \basel{\hat{f}}{i}(\xx)\vert < \delta$,
 for $1 \leq i \leq m$, for every $n \geq \nu$. Thus, every weak basic open set
 centered  at $\xx$ intersects ${\mathcal M }$, and $\xx$ is a limit point of
 ${\mathcal M}$ with respect to the weak topology.                     \qed
\\

\begin{theorem}
{\bf{(Weak Denseness of a Countable Set)}}\tab   
\label{Weak-Denseness-of-a-Countable-Set}
 Let ${\mathcal B}$ be a Banach space with a countable basis
 $\{\basel{\xi}{i} \,:\, i\in \naturalNums\}$, and
 $\basel{\mathcal M}{\mathbb Q} =  \bigcup_{n \in \naturalNums}
  \basel{\mathcal L}{\mathbb Q} (\{1,\, ... , \,n\})$,
 where $\basel{\mathcal L}{\mathbb Q} (\{1,\, ... , \,n\})$
 is the countable subset obtained by collecting linear combinations of
 vectors in $\{\basel{\xi}{i} \,:\, 1 \leq  i \leq  n\}$, with rational
 number  coordinates,  for  $n \in \naturalNums$.
 Then, $\basel{\mathcal M}{\mathbb Q}$ is weakly dense in ${\mathcal B}$.
 \end{theorem}
\proof $\basel{\mathcal M}{\mathbb Q}$ is strongly dense
 in ${\mathcal M}$, where ${\mathcal M}$ is as defined in
 Theorem \ref{Weak-Denseness-of-Finite-Linear-Combinations}
and  weakly dense in ${\mathcal B}$.   \qed
\\

\begin{Corollary}
\label{Corollary-to-Weak-Denseness-of-a-Countable-Set}
Let ${\mathcal B}$, ${\mathcal M}$ and $\basel{\mathcal M}{\mathbb Q}$
be as in Theorems \ref{Weak-Denseness-of-Finite-Linear-Combinations}
and \ref{Weak-Denseness-of-a-Countable-Set}. For any continuous
linear functional $\hat{g}$ defined on
 $\basel{\mathcal M}{\mathbb Q}$ , such that
$\vert \hat{g}(\xx) \vert \leq \norm{\xx}{\mathcal B}$ ,
for $\xx \in  \basel{\mathcal M}{\mathbb Q}$ ,
there is a unique continuous extension $\hat{h}$
of $\hat{g}$ to ${\mathcal M}$, and $\hat{h}$ satisfies
 $\vert \hat{h}(\xx) \vert \leq \norm{\xx}{\mathcal B}$ ,
 for every $\xx \in  {\mathcal M}$.
\end{Corollary}
\proof $\basel{\mathcal M}{\mathbb Q}$ is strongly dense
 in ${\mathcal M}$.  \qed
 \\
 
 For a property that depends continuously with respect to the weak topology,
 in order to show that the property holds for all of  ${\mathcal B}$, 
 it suffices to show that the same holds for ${\mathcal M}$.
\\

\textcolor{titlecolor}{
\section{Conclusions} 
}
The projection maps and dual basis linear functionals a Banach space 
with a countable vector space basis are shown to be continuous open maps.
The dual basis linear functionals form a basis for the dual spacethe dual space,
and the double dual basis linear functionals form a basis for the double dual space.
A specific criterion for reflexivity is presented.

 \textcolor{refscolor}{
}
\end{document}